\documentclass[11pt]{amsart}

\usepackage{rotating}
\usepackage{amssymb}
\usepackage{amsxtra}
\usepackage{graphics}
\usepackage{latexsym}
\usepackage{xcolor}
\usepackage{amsmath}
\usepackage{amssymb,amsthm,amsfonts}
\usepackage{mathtools}
\usepackage{amscd}
\usepackage[arrow, matrix, curve]{xy}
\usepackage{syntonly}
\usepackage{multirow}
\usepackage{fancyhdr}
\usepackage{hyperref}

\usepackage{etex}
\usepackage{tikz-cd}
\usetikzlibrary{matrix, calc, arrows}
\usepackage{braket}
\usepackage{adjustbox}
\usepackage[official]{eurosym}

\theoremstyle{plain}

\theoremstyle{definition}

\numberwithin{equation}{thm}

\newcommand{\rk}{{\rm rank}}

\newcommand{\en}{End}

\usepackage{geometry}

\begin{document}\textcolor{white}{.}\\

\centerline{\bf Higgs Bundles in Geometry and Arithmetic}\vspace{.2cm}
\centerline{\bf Kang Zuo}
\vspace{.2cm}\textcolor{white}{.}\\	
This is an expository note based on my lecture in 2020 ICCM annual meeting at the
 University of Science and Technology
of China in Hefei.\\[.2cm]
A Higgs bundle on a complex manifold $X$ is a pair $(E,\theta)$ consisting of a holomorphic vector bundle $E$ on $X$ and an $\en(E)$-valued 1-form $\theta$ satisfying the condition $\theta\wedge\theta=0$. The form $\theta$ is called the \emph{Higgs field}. A major development in the field was made by N. Hitchin \cite{H} and C. Simpson \cite{S1} in the so-called Hitchin--Simpson correspondence, a powerful tool in complex algebraic and analytic geometry. They followed the physicists in solving the Yang--Mills--Higgs partial differential equations on semistable Higgs bundles. Generalizations of this approach to more abstract settings in algebraic geometry have been highly successful \cite{LSYZ, EG, Lan, Ar, KYZ1, KYZ2, YZ}. \\[2mm]
In the project proposed here we introduce the notions of  deformation Higgs bundle and  Riemann--Finsler metric on the moduli space of polarized varieties. We also use the Higgs--de Rham flow in the $p$-adic setting. These are the key novelties in our program. These tools enable us to attack the following problems:\\[.2cm] $\bullet$ {\sl The Shafarevich conjecture on the finiteness of isomorphism classes of families of higher dimensional varieties.} \\[.2cm] $\bullet$ {\sl A folklore conjecture on the bigness of the fundamental group of moduli spaces of smooth projective  varieties with semi-ample canonical line bundles.}\\[.2cm] These conjectures are major current roadblocks preventing the progress in the Shafarevich program.\\[.2cm] In the $p$-adic nonabelian Hodge Theory, we develop and explore further a theory of Higgs bundles on varieties over $p$-adic fields. Two directions of applications are\\[.2cm] $\bullet$ {\sl The Faltings' $p$-adic Simpson correspondence.} \\[.2cm] $\bullet$  {\sl The construction of motivic local systems over $p$-adic curves in connection with Drinfeld's work on the Langlands program.} \\[.2cm]Below, we explain the above fundamental notions in our approach and the problems we will consider. \\[.2cm]
{\bf I.  Higgs Bundles  on Moduli spaces   of Manifolds and  the Shafarevich Program.}\\
The Shafarevich conjecture over function fields  of one variable,  proved by Parshin and Arakelov \cite{A}, states that: {\bf a.} If a smooth algebraic curve $U$  parametrizes a non-isotrivial  family of smooth projective curves, then $U$ is a hyperbolic curve, i.e.  $U$  carries a metric of constant negative curvature.
{\bf b.} There are only finitely many   non-isotrivial families parametrized by $U$ up to isomorphism.\\[.2cm]
In algebraic geometry, moduli spaces are spaces whose points parametrize solutions to a specific geometric problem. One of the most important moduli spaces is  the   moduli space $M_h$ parameterizing algebraic manifolds with semi-ample canonical line bundle and with  a fixed Hilbert polynomial $h$.  A family  $f\colon V\to U$ of $n$-dimensional complex manifolds gives rise to a classifying map  $\phi\colon U\to  M_h$ and the family is called of \emph{maximal variation} if $\dim U=\dim\phi (U).$ Kodaira and Spencer studied the variation of the complex structure of the family $f$ by introducing   the  notion of the Kodaira--Spencer map
$\tau^{n,0}\colon  T_U\to R^1f_*T_{V/U}$. This map can be extended to a  log-smooth  compactification $f\colon X\to Y$, i.e., to the logarithmic Kodaira-Spencer map $\overline{\tau}^{n,0}:T_Y(-\log S)\to R^1f_*T_{X/Y}(-\log\Delta).$ 	
Griffiths  \cite{G1,G2,G3}  introduced the notion of polarized variation of Hodge structures on  $U$ by  examining how the Hodge structures on the Betti cohomology of the fibers of $f$ vary.  Simpson  \cite{S1}  introduced  a basic notion in nonabelian Hodge theory, the so-called  system of Hodge bundles  $(E,\theta)$, to be the associated graded Higgs bundle of a  polarized variation of Hodge structures. These theories are of fundamental importance  in the study of the  global geometry of  families. The works of Viehweg--Zuo \cite{VZ-0,VZ-2,VZ-1,VZ-3}, combined the above-mentioned theories by  extending the logarithmic Kodaira-Spencer map  to the  higher direct images of wedge products $T^q_{X/Y}(-\log \Delta):=\wedge^q T_{X/Y}(-\log \Delta)$ of the relative log  tangent sheaf $\tau^{n-q,q}:T_Y(-\log S)\otimes R^{q}f_*T^q_{X/Y}(-\log\Delta)\xrightarrow{}R^{q+1}f_*T^{q+1}_{X/Y}(-\log\Delta)$ and constructed  the \emph{logarithmic (graded) deformation Higgs bundle}
$(F,\tau):=( \bigoplus_q R^{q}f_*T^q_{X/Y}(-\log \Delta), \, \bigoplus_q\tau^{n-q,q})$
for the family $f\colon  X\to Y$.  Moreover, Viehweg--Zuo  constructed  a comparison map 
$\rho\colon (F,\tau)\to (E,\theta)\otimes A^{-1},$ 
where $(E,\theta)$ is the Higgs bundle associated to the  variation of Hodge structures on the middle cohomology  of  a new family  $g\colon Z\to Y$  built from  a cyclic covering of $X$ by taking roots out of sections of the relative pluricanonical linear system twisted by an anti-ample line  bundle $A^{-1}$ on $Y$.  
Via the  maximal non-zero iteration $ \tau^m: S^mT_Y(-\log S)\xrightarrow{}
\text{ker}(\theta^{n-m-1,m+1}) \otimes A^{-1}$ of the Kodaira--Spencer map of $(F,\tau)$ (also known as the Griffiths--Yukawa coupling), where $\text{ker}(\theta^{p,q})$ is the kernel of the Higgs field $\theta^{p, q}$ which is semi-negative \cite{Z},
one obtains a big subsheaf 
$\mathcal A\subset S^m\Omega^1_Y(\log S).$
In the analytic setting, one defines   a complex Finsler pseudometric $ds^2_{\mathbb{ C} \,F}$ on $T_Y(-\log S)$ by taking the $m$-th root out of the product of the Hodge metric and the Fubini-Study metric  on $ \text{ker}(\theta^{n-m-1,m+1}) \otimes A^{-1}$   via the map $\tau^m$.  One shows that  
$ds^2_{\mathbb{ C} \,F}$  has the holomorphic sectional curvature  bounded above by a negative constant.
The negative holomorphic sectional curvature derived as above plays the same role as the negative holomorphic sectional curvature associated to horizontal period maps in Hodge theory.  It  implies Brody and Kobayashi hyperbolities of $U$ in \cite{VZ-1}, \cite{PTW} and \cite{De}. Also,   the big Picard theorem was recently proven  in \cite{DLSZ}. \\[.2cm]
{\bf I.1.Finiteness  of Families and Geometric Characterization of Non-Rigid Families.} The \emph{Finiteness problem} proposed by the Shafarevich program for the isomorphism classes of families of higher dimensional varieties may be reduced to two subproblems: the boundedness of  families and the rigidity of  families. The boundedness of  families over a fixed base has been proven by various people,  see for example \cite{F, VZ-0,  KL} . The \emph{ rigidity problem} is crucial for solving the conjecture on the finiteness and  is subtle. Recall that a family $f:V\to U$ is  called non-rigid if the family can  be non-trivially extended over $U\times \Delta$, where $\Delta$ is a disc. One can construct non-rigid families of higher dimensional varieties by simply taking the product of two families with lower-dimensional fibers. 
Under certain additional assumptions we expect the rigidity holds true.  As an initial step, a recent work of Javanpeykar--Sun--Zuo  \cite{JSZ} shows that the  $N$-pointed 
Shafarevich conjecture is true: let $C$ be a smooth projective curve of genus  $g$ and $S$ a divisor on $C$ with given points $\{c_1,\dots, c_N\}\subset C\setminus S $.  Then there are only finitely many isomorphism classes of families $f:X\to C$ of $n$-folds with semi-ample canonical line bundles and bad reduction at $S$ and with the fixed isomorphism classes of the 
fibers over $\{c_1,\dots, c_N\}$  for 
$N\geq {n-1 \over 2}(2g-2+\#(S))$. Currently,  Javanpeykar-Lu-Sun-Zuo are making new progress  on this type of problems and we expect that the one-pointed Shafarevich conjecture holds true. The proofs crucially rely on the deformation Higgs bundle and the complex Finsler metric whose curvature is bounded above by a negative constant.
At the final  step  we  study  the rigidity problem in the original Shafarevich program.  
Given a non-rigid family $ f \colon X\to C$  over a fixed base curve $C$ and with a fixed degeneration locus $S_C$,  the boundedness of families over a fixed base with a fixed degeneration locus  implies that $ f $ extends to a non-trivial family  over a product space  $f\colon X\to C\times T$, where $T$ is a projective curve and  with the degeneration locus
$S_C\times T\cup C\times S_T.$\\[.15cm]
{\bf Problem 1.} {\bf a.}  Show that  the deformation Higgs bundle for the family $f\colon  X\to C\times T$ with the projections $p\colon C\times T\to C$ and $ q\colon  C\times T\to T$  decompose  in the following form
$ (F,\tau)=p^*(F, \tau)_C\otimes q^*(F,\tau)_T\oplus (F,\tau)'',$
where $(F, \tau)_C$ and $ (F,\tau)_T$ are graded logarithmic Higgs bundles over $(C, S_C)$ and $(T, S_T).$  \\
{\bf b.}	Show that a family $f \colon X \to C\times T$  of canonically polarized manifolds  over a product base with maximal variation  decomposes 
as the product  of  two  families  $ f_{C}  \colon X_{C}\to C$  and  $ f_{T} \colon X_{T}\to T$ with lower dimensional fibers after taking a ramified cover of $X$ (up to blow-ups and blow-downs). \\[.15cm]
The statement in  Problem 1.b.  should hold  in general only for fibers of general type as Faltings has found   examples of  non-rigid families of abelian varieties whose  generic fibers are simple  abelian varieties. \\
The investigation  of non-rigid families of hypersurfaces in projective spaces  will be the first non-trivial test for Problem 1. In \cite{VZ-3}, Viehweg--Zuo have constructed non-rigid families of hypersurfaces in $\mathbb P^n$. For simplicity, here  we consider hypersurfaces in $\mathbb P^3$  of dgree $d \geq 4$.
Any set of $d$ (pairwise) distinct points in $\mathbb P^1$ corresponds to a homogeneous polynomial in two variables of degree $d$ which has this set as its set of roots. For any two such homogeneous polynomials $f_d(x_1,x_2)$ and $ g_d(y_1,y_2)$, the equation $ f_d(x_1, x_2)+g_d(y_1, y_2)=0$ defines a smooth hypersurface $X_{f_d, g_d}$ in $\mathbb P^3$. By varying the roots of $f_d$ and $g_d$ (and thus varying $f_d$ and $g_d$), one obtains a family $ f: X\to \Sigma_d\times \Sigma_d$ in $\mathbb P^3$, where  $\Sigma_d$ is the moduli space of $d$ distinct points in $\mathbb P^1.$ Viehweg-Zuo have observed that the family $f$ is, up to ramified covers and blow-ups and blow-downs, the product of the family  $C_{f_d}\to \Sigma_d$  defined by the equation $ x_3^d=f_d(x_1,x_2)$  and the family $C_{g_d}\to \Sigma_d$ defined by $ y_3^d=g_d(y_1,y_2)$. In fact we have enough evidence to be convinced that
\\[.2cm]
{\bf Conjecture 2.}  Any smooth  fiber of  a  non-rigid family of hypersurfaces in $\mathbb P^3$  of degree $d$ is defined
by an equation   of the form $f_d(x_1,x_2)+g_d(y_1,y_2)=0$ up to a projective transformation.\\[.2cm]
Viehweg--Zuo have found a close relationship between the relative graded Jacobian ring and the system of Hodge bundles by showing that  the  graded pieces of Jocobian ring of a hypersurface $X\subset \mathbb P^n$  coincide with  the eigenspaces of the $\mathbb Z/d\mathbb Z$-action on the middle-dimensional Hodge cohomology of the cyclic cover 
of $\mathbb P^n$ ramified over $X$ as modules together with the Kodaira-Spencer multiplication. Using Deligne's result on tensor decomposition of a polarized variation of Hodge structures over a product base one  shows that  a decomposition of the base of a family $f\colon X\to C\times T$  of hypersurfaces in the projective space $\mathbb P^n$  leads to a tensor product decompostion of the reletive Jacobian ring as modules together with the Kodaira-Spencer multiplication. Moreover,   Problem 1.a. has a positive answer  in this case.
It is well-known that the Jacobian ring of  a smooth hypersurface determines the  isomorphism type of the  hypersurface itself. We hope that the above decomposition of the relative Jacobian ring (as a module) as well as the decomposition of the deformation Higgs bundle will lead to some significant geometric consequences on the fibers, as in Conjecture 2. Motivated by   Problem 1.b. we pose the following problem which gives a criterion for rigidity.\\[.15cm]
{\bf Problem  3.}  Show that  any family $f\colon V\to U$  of maximal variation  and containing a fiber  with the generically ample cotangent sheaf  is   rigid.\\[.15cm]
{\bf  I.2. Bigness of Fundamental Groups of Moduli Spaces of Projective Manifolds.}   Milnor \cite{Mil68} introduced  a growth function $\ell$ associated to a finitely generated group $G$ as follows: for each positive integer $s$ let $\ell(s)$ be the number of distinct group elements which can be expressed as words of length $\leq s$ with a fixed choice of generators and their inverses. 
Milnor proved that $\pi_1(M)$ of a compact Riemannian manifold $M$ with all Riemannian sectional curvatures less than zero has exponential growth, i.e.,  $\ell(s) \geq a^s$  for some $a>1.$ The proof relies on  G\"unther's volume comparison theorem on the exponential growth of the volume of the geodesic ball  on the universal cover $\tilde M$   \cite{Gue}.  We  propose the following problem. \\[.15cm]
{\bf Problem 4.}  Let $U$ be a base space parameterizing polarized manifolds with semi-ample line bundle and with maximal variation.  Prove that  $\pi_1(U)$  grows at least exponentially.   \\[.15cm]
Similar to the approach of proving the complex hyperbolicity  on $U$ we  would like to  construct a Riemann-Finsler  metric  on $U$  via  the maximal non-zero iteration of Kodaira--Spencer map, whose  curvature has certain negativity properties.
Wu--Xin  generalized Milnor's theorem to the Finsler setting by showing that the G\"unther's volume comparison theorem  holds for Riemann-Finsler metric  with negative flag curvature  (an analogue of the Riemannian sectional curvature in Finsler geometry) \cite{CZ, WX}. In our situation the complex Finsler metric (the sum of  several such metrices from different cyclic covers)  naturally induces a non-degenerated  Riemann-Finsler metric $ds^2_{\mathbb{R}\text{Fin}}$ on $U$.  We are aware that  in general the Riemannian curvature decreasing principle  does not hold true for real submanifolds.  Very recently, together with S. Lu and R.R.  Sun, we  observed that the fact that pluriharmonicity  of the  composition of the horizontal period map   with the projection to the symmetric space of non-compact type  implies decreasing Riemannian  curvature still holds true in a weak form and we expect that this weak form of the negative sectional curvature property can show G\"unther's volume comparison inequality and  get the solution of   Problem 4 by Milnor's original argument.
\\[.4cm]
{\bf II. Higgs Bundles in $p$-adic Nonabelian Hodge Theory and Motivic Local Systems.}\\
In the $p$-adic world, Fontaine  and Faltings  introduced the notion of  Fontaine--Faltings module
${\mathcal MF}^{\nabla}({\mathcal X}):=\{ (V,Fil^\bullet,\nabla,\Phi)\}$  over a
smooth  proper scheme $\mathcal {X}$ over a $p$-adic number ring $\mathcal{O}_{K_p}$ (more precisely on the $p$-adic completion  $\hat {\mathcal X}$). A Fontaine--Faltings module consists of a filtered de Rham bundle  endowed with a relative Frobenius  $\Phi$.  The  Fontaine-Laffaille-Faltings $p$-adic  Riemann-Hilbert
correspondence produces a functor  $\mathbb{D} $ sending a Fontaine--Faltings module to a crystalline representation of  $\pi^\text{\'et}_1(X_K)$  (\cite{FL},\cite{F3}).
Faltings (\cite{F4})  has proposed a  $p$-adic Simpson correspondence  to describe the so-called generalized representations of the geometric fundamental group $\pi^\text{\'et}_1(X_{\bar{\mathbb Q}_p})$ in terms of Higgs bundles.  Faltings asked whether semi-stable Higgs bundles  $(E,\theta)$  with trivial Chern classes $c_i(E)=0$  correspond precisely  to genuine representations. Scholze \cite{Sch} has made exciting progress in $p$-adic Hodge theory. He introduced  the pro\'etale site  $X_\text{pro\'et}$ of $X$ as a refinement of the usual \'etale topology on $X$; using this topology, he defined de Rham representations as being associated to an $\mathcal {O}\mathbb{ B}_\text{dR}$-module equipped with a Hodge filtration. By taking the category  $\mathcal {HIG} $ of semi-stable graded Higgs bundles over  ${\mathcal X},$  Lan--Sheng--Zuo \cite{LSZ1} introduced the notions {\sl Higgs--de Rham flow} and {\sl periodic Higgs bundle}, a  $p$-adic  analogue of Yang--Mills--Higgs equation over the complex numbers. They established  a $p$-adic Simpson correspondence between  the category of  stable periodic Higgs bundles  $ \mathcal {HIG}^{\text{per}}  \subset  \mathcal {HIG}  $ and the category  $\mathcal {REP}^\text{cry}$ of crystalline representations of $\pi^\text{\'et}_1(X_K)$ when $X$ has good reduction.
The construction
relies on the fundamental work of Ogus--Vologodsky \cite{OV} on Cartier transform in characteristic $p$,  and the work of Simpson on the existence of grading semi-stable Hodge filtration \cite{S2}.  
The notion of Higgs--de Rham flow has already  applications in both algebraic and arithmetic geometry, see recent works by \cite{LSYZ, EG, Lan, Ar, KYZ1, KYZ2}.\\[.2cm]
{\bf II.1. $p$-adic Simpson Correspondence and Semi-Stable Higgs Bundles.}  
For simplicity we just work with a smooth projective curve $\mathcal X$ over $W(k)$. Given a Higgs bundle $(E,\theta)$ over  $\mathcal{X}$  with nilpotent Higgs field of degree $<p$. By Faltings' Simpson correspondence for integral generalized representation of $\pi_1^\text{geo}$. We believe that the heart part for proving Faltings conjecture 
lies in the characteristic $p$ case, and the stability of Higgs bundle plays the crucial role.\\[.2cm]
{\bf II.1.1.  Faltings' Conjecture over  Characteristic $p$ and Preperiodic  Higgs de Rham Flow.}\\
Given  a  semi-stable 
Higgs bundle $(\bar E,\bar \theta)_0$  on  $X_k=\mathcal{X}$ (mod $p$) of degree zero and with a  nilpotent Higgs field of degree $<p$, we obtain a preperiodic Higg-de Rham flow by runing semi-stable  Higgs de Rham flow over $X_k$  with the initial  $(\bar E,\bar\theta)_0$ 
{\tiny 	
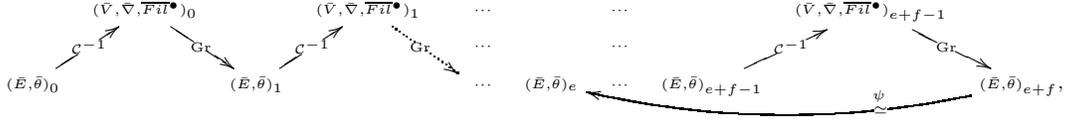
\begin{equation*}
	\xymatrix@C=2mm@R=5mm{ 
		\,
		& {\scriptstyle (\bar V, \bar \nabla, \overline{Fil}^\bullet)_0  } \ar[rd]|{ \text{Gr} }   
		& 
		& {\scriptstyle (\bar V, \bar \nabla, \overline{Fil}^\bullet)_1  } \ar@{.>}[rd]|{ \text{Gr} }   
		& {\scriptstyle \quad \cdots } \ar@{}[d]|{\quad\cdots}  
		&
		& {\scriptstyle \cdots  }  \ar@{}[d]|{\cdots}  
		&
		& {\scriptstyle (\bar V, \bar \nabla, \overline{Fil}^\bullet)_{e+f-1}   }  \ar[rd]|{ \text{Gr} }   
		&
		\\
		{\scriptstyle (\bar E, \bar \theta)_0   } \ar[ru]|{ \mathcal C^{-1} }   
		&
		& {\scriptstyle (\bar E, \bar \theta)_1 }   \ar[ru]|{ \mathcal C^{-1} }          
		&
		& {\scriptstyle \quad \cdots     }
		& {\scriptstyle  (\bar E, \bar \theta)_e    }
		& {\scriptstyle \cdots   }
		& {\scriptstyle (\bar E, \bar \theta)_{e+f-1}  }  \ar[ru]|{ \mathcal C^{-1} }
		& 
		& {\scriptstyle  (\bar E, \bar \theta)_{e+f}}, \ar@/^1pc/[llll]|(0.3){ \simeq}_(0.3){ \psi} 
	}
	\end{equation*} }
where the Hodge filtrations in the  semi-stable  de Rham bundles  are  grading semi-stable.  The Lan-Sheng-Zuo's  functor   $\alpha$ over $k$ is defined by  sending the direct sum of the periodic part  of the flow with the $f$-periodic isomorphism $\Psi$ 
to the direct sum of their inverse Cartier transform endowed with the Hodge filtration
$(\bigoplus_{i=0}^{f-1}(\bar E,\bar \theta)_{e+i}), \Psi)\overset{\alpha}\mapsto(\bigoplus_{i=0}^{f-1}(\bar V,\bar \nabla, \overline{Fil}^\bullet)_{e+i}, \mathcal{C}^{-1}(\Psi)),$
which  is  indeed a $p-$torsion Fontaine-Faltings module, endowed with a natural  $\mathbb{ F}_{p^f}$-action, and therefore  corresponds to a $\mathbb F_{p^f}$-crystalline representation  $\mathbb L^\text{cry}_{(\bar E,\bar \theta)_e}.$ This correspondence restricted to the geometric $\pi_1$ coincides with Faltings' Simpson correspondence.
For a morphism $\sigma: \mathcal Y\to \mathcal X$ and a small Higgs bundle $(E,\theta)/\mathcal X$, Faltings introduced the  notion of  the twisted pull-back $\sigma^\circ (E,\theta)$ of $(E,\theta)$ in \cite{F4}, which is compatible with the pull-back of generalized representations under Faltings' Simpson corrspondence.\\[.2cm]
{\bf Conjecture 5.}    Given a preperiodic  Higgs-de Rham flow as above, there  exists a   morphism $\sigma: \mathcal Y\to\mathcal X$ such that $\sigma^\circ (\bar E,\bar\theta)_0\simeq \sigma^\circ(\bar E,\bar\theta)_e.$
\\[.2cm]	
{\bf Remark 6.}  Conjecture 5 implies that a semi-stable nilpotent Higgs bundle 
over $X_k$ of degree zero corresponds to an $\mathbb{F}_{p^f}[\pi]/(\pi^{p^{e-1}(p-1)})$-geometric local system.
\\[.2cm] 
{\bf II.1.2. Lifting of Higgs Bundles and Local Systems at Truncated Level.}
Given a Higgs bundle $(E,\theta)$ over $\mathcal X_{\mathcal O_{\mathbb C_p}}$ such that
$(\bar E,\bar\theta)$ over $X_{\mathcal O_{\mathbb C_p}/p\mathcal O_{\mathbb C_p}  }$ corresponds to a genuine geometric local system $\bar {\mathbb L}$, by taking a Galois base change  $\sigma: \mathcal Y_{\mathcal O_{\mathbb C_p}}\to \mathcal X_{\mathcal O_{\mathbb C_p}}$, \'etale on the generic fiber killing $\bar {\mathbb L},$  we may assume Faltings' twisted pull-back $\sigma^\circ(E,\theta)$ on $Y_{\mathcal O_{\mathbb C_p}/p\mathcal O_{\mathbb C_p}  }$ is isomorphic to  the trivial Higgs bundle $(\mathcal O^r,0)$.  One can measure   the difference between
$\sigma^\circ(E,\theta)$ and  $(\mathcal{O}^r, 0)$ over 
$Y_{ \mathcal{O}_{\mathbb C_p}/p^2\mathcal{O}_{\mathbb C_p}}$ by a class $\eta_\text {Hig} $ in $ H^1_\text{Hig}(Y_ {\mathcal{O}_{\mathbb C_p}/p\mathcal{O}_{\mathbb C_p}}, \mathcal{E}nd(\mathcal O^r, 0))$ as the torsor space for lifting  Higgs bundles 
from modulo  $p$ to modulo  $p^2$  over the trivial Higgs bundle.\\[.2cm]
{\bf Problem 7.} Work out the Faltings'  $p$-adic Simpson correspondence between the torsor space  for lifting  Higgs bundles and the torsor space for lifting geometric local systems and show that $(E,\theta)$ corresponds to a genuine geometric local system $\mathbb L.$\\[.2cm]
Given a generic punctured hyperbolic curve $(C,S)/\mathbb F_q,$ in \cite{LSYZ} we showed  that there exists a so-called {\sl canonically lifted} $(\mathcal C,\mathcal S)/W(\mathbb F_q)$ regarding that  the uniformizing Higgs bundle $\theta: L\xrightarrow{\simeq} L^{-1}\otimes \Omega^1_{\mathcal C}(\log\mathcal S)$ corresponds to a crystalline local system on $(C\setminus S)/W(\mathbb F_q)[1/p],$
where $L$ is the logarithmic Theta characteristic of $(\mathcal C,\mathcal S)$. This theorem shall be considered  as the Higgs bundle incarnation
of   Mochizuki's $p$-adic Teichm\"uller theory \cite{Mo}.
The techniques developed in attacking problems  in II. 1.1. and  II.1.2  shall allow us  to work out  the far-reaching generalization  for $p$-adic Teichm\"uller theory for arbitrary 
$p$-adic  punctured hyperbolic curves.\\[.15cm]
{\bf Problem 8.} ($p$-adic Teichm\"uller theory) Show that  the  uniformizing Higgs bundle over a $p$-adic punctured hyperbolic curve $(C,S)/K_p$ corresponds to a geometric $\mathbb C_p$-local system $\rho^\text{uni}$, which is invariant under the natural $\text{Gal}(\bar{K}_p/K_p)$-action.  Furthermore, if  $(C,S)\not\simeq (C', S')$ then $\rho^\text{uni}\not\simeq \rho^{'\text{uni}}.$
\\[.15cm]
{\bf II.1.3.  Galois-Action on Geometric $\pi_1$ and a $p$-adic Analogue of $\mathbb C^*$-action on Higgs Fields.} The action of Galois group $\text{Gal}(\bar K_p/K_p)$ on $X_{\bar K_p}$ induces a natural action on the category of generalized representations. In \cite{YZ} by carefully checking the construction of Faltings' Simpson correspondence one finds that the corresponding action on the category of Higgs bundles can be basically described as the  Galois action on the usual category of  Higgs bundles   with the  extra action on
Higgs fields  twisted by a 1-cocycle induced by the element
$\xi$ in the  Fontaine periods ring $ B^+_\text{dR}.$ Consequently, we show that a generalized representation corresponding to a graded Higgs bundle over $K_p$ is $\text{Gal}(\bar K_p/K_p)$-invariant and conversely, a Higgs bundle corresponding to a Galois-invariant generalized representation is  precisely graded and defined over  a finite extension of $K_p$
for $\rk E\leq 2$, and with a nilpotent Higgs field  in the general case. Together with the expected results from II.1.1. and II.1.2. we shall get  a better understanding of the Galois action on the geometric $\pi_1$ and, in particular,  Grothendieck's anabelian goemetry.\\[.2cm]
{\bf II.2. Rank-2  $p$-adic   Motivic Local Systems on a Punctured Curve  and Langlands Correspondence  over Function Fields.} 
Let $(C, S)$ be a smooth projective curve over  complex numbers  together with a set of $n$ punctures $S$. Simpson shows that
the category of abelian schemes over $C$ of $\text{GL}_2$ type  with bad reductions over $S$ is equivalent to the category of rank-2 motivic Higgs bundles on $C$ with logarithmic-parabolic structure on $S$.  We plan to study rank-2 motivic Higgs bundles
using $p$-adic nonabelian Hodge theory and Langlands correspondence over function fileds in char. $p.$
Let  $(\mathcal C, \mathcal S)_{W(k)}$   denote  a smooth projective curve  with a set of $n$ punctures $\mathcal S$
over $W(k)$,  and   $\mathcal {HIG} $ denotes the moduli space of rank-2  graded stable Higgs bundles  over $(\mathcal C,\mathcal S)$ with log-parabolic structure over $S$.
Clearly a motivic Higgs bundle  on an arithmetic scheme must be  periodic over all  unramifiled places.
Conversely, we ask: \\[.2cm]
{\bf Problem 9.}  Study  those rank-2  log periodic Higgs bundles over
$(\mathcal C,\mathcal S)/W(\mathbb F_q)$, which are  motivic.
\\[.2cm]
In \cite{SYZ} one considers 
rank-2 stable  logarithmic Higgs bundles over the  projective line $(\mathcal P^1,\mathcal S)$ 
with $n$ punctures
and with the parabolic structure attached to one puncture in  $\mathcal S$ of weight $\{1/2,\, 1/2\}$.  The moduli space $\mathcal {HIG}$ contains irreducible components of dimensions $n-3,\,  n-5,\, n-7, \cdots$.  Sun--Yang--Zuo  constructed infinitely many $\text{GL}_2(\mathbb{Z}_p^\text{ur})$-crystalline local systems  over   the generic fibre $( \mathbb{P}^1\setminus S)_{\hat {\mathbb Q}_p^\text{ur}}$  by running Higgs--de Rham flow on  $\mathcal {HIG} $  and showing that the set of periodic Higgs bundles over $\mathbb Z_p^\text{ur}$ is Zariski dense in the components of $\mathcal {HIG}$  of maximal dimension. That should be also true for a general punctured hyperbolic curve and  follow from  the deformation theory developed in \cite{KYZ3}. We raise:
\\[.15cm]
{\bf Problem 10.}  
Given a punctured hyperbolic curve $(\mathcal C,\mathcal S)/_{W(\mathbb F_q)},$ show that  the set of periodic Higgs bundles
over $W(\mathbb F_{p^n}),\, n\in \mathbb N$  is Zariski dense
in $\mathcal {HIG} $.   Are they motivic?\\[.15cm]
Problem 10 is a $p$-adic analogue of  a very recent program of Esnault-Kerz \cite{EK}.
They conjecture that motivic points are Zariski dense in the
character variety of $\ell$-adic  local systems.  \\
For $\mathcal S=\{0,1,\infty,\lambda\}$,  the moduli space  $\mathcal {HIG}$
parameterizes rank-2 graded stable  Higgs bundles over $(\mathcal  P^1,\mathcal S)/W(\bar k)$ of the form 
$(\mathcal O\oplus\mathcal O(-1),\,  \mathcal O\xrightarrow{\theta\not=0}\mathcal{O}(-1)\otimes \Omega^1_{\mathcal P^1}(\log S))$
and  $\theta$ has a   single zero $(\theta)_0\in \mathcal{P}^1$. One makes an identification 
$\mathcal {HIG} = \mathcal P^1$
by sending $(E,\theta)\in \mathcal {HIG} $ to
$(\theta)_0\in \mathcal P^1$.  We take  the elliptic curve $(\mathcal E_\lambda\colon y^2 = x(x-1)(x-\lambda),\,0:=\infty)$ as the double cover of $\mathcal P^1$ ramified over $\mathcal S.$  Recently Krishnamoorthy, Yang and Zuo are making  new progress. We show  that Problem 10  has a positive answer  for $(\mathcal P^1,\{0,1,\infty,\lambda\})$ for $p<50$ and $\mathcal E_\lambda$ is a supersingular elliptic curve.\\[.2cm]
{\bf Conjecture 11.} ( \cite{SYZ})  The set  of periodic Higgs bundles over $(\mathcal P^1,\{0,1,\infty,\lambda\})/W(\mathbb F_q),\, n\in \mathbb N$ carries a "group law"  from the elliptic curve $\mathcal E_\lambda, $ i.e. a Higgs bundle  $(E,\theta)\in \mathcal {HIG}(W(\mathbb F_{p^n}))$ is $f$-periodic 
if and only if  the preimage of the zero of the Higgs field $\pi^{-1}(\theta)_0\in \mathcal E_\lambda$ is a $(p^f\pm 1)$-torsion point. \\[.15cm]
{\bf Langlands Correspondence over Function Fields of Characteristic  $p$.}
Given a  punctured  curve $(\mathcal C, \mathcal S)_{W(\mathbb{F}_q)}$, and picking up a field isomorphism $\sigma\colon \overline{\mathbb Q}_p\rightarrow \overline{\mathbb Q}_l$, by Deligne's conjecture on $p$-to-$l$ companions solved by  Abe \cite{Ab}  and via overconvergent $F$-isocrystals one obtains an inclusion  from the set of (modulo $p$) stable logarithmic  periodic Higgs bundles $ \mathcal {HIG}^\text{per} \hookrightarrow \mathcal{REP}^{\ell\text{-adic}},$ the set  of  irreducible $\text{GL}_2(\bar {\mathbb Q}_\ell)$ -local systems over $ ( C \setminus S)_{\mathbb F_q}$.  By Drinfeld \cite{D2} the $\ell$-adic local system correspoinding to $(E,\theta)$ together with its full companions will be realized by  an abelian scheme over $(C\setminus S)_{\mathbb F_q}$ of $\text{GL}_2$-type. One studies the existence of grading stable Hodge filtrations attached to realizations of  $F$-crystals  over $(\mathcal C,\mathcal S)_{W(\mathbb F_q)},$ applies  the logarithmic Grothendieck-Messing deformation theorem and shows that this abelian scheme lifts over $W(\mathbb F_{p^n})$ for the case $(\mathcal P^1,\{0,1,\infty, \lambda\}).$ In general, we pose 
\\[.2cm]
{\bf Problem 12.} Find  condtions  for abelian schemes of $\text{GL}_2$-type over $(C,S)_{\mathbb F_q}$ arising from Langlands correspondence  for rank-2 $\ell$-adic local systems to be liftable over $W(\mathbb F_q)$.
\\[.1cm]
{\bf $p$-adic Nonabelian Hodge Theory and Dynamical Systems over $\mathbb F_p$ and Mixed Characteristic.}
Conjecture 11  is   motivated by 
Sun-Yang-Zuo's work \cite{SYZ}.   The Higgs-de Rham flow on $(\mathbb P^1, \{0,1,\infty,\lambda\})$  induces a self map $\psi_{\lambda}$ on the moduli space $\mathcal{HIG}\,(\simeq \mathbb P_{\mathbb F_q}^1$) of rank-2 logarithmic  Higgs bundles  with 
the parabolic structure of the weight (1/2,\,1/2) attched to a puncture. As is explained in {\bf II.1.1}, a point Higgs bundle $z\in \mathcal{HIG}$
corresponds to a crystalline local system on $\mathbb P^1\setminus\{0,1,\infty,\lambda\}$ iff $z$ is a periodic point of $\psi_{\lambda}.$\\[.15cm]
{\bf Conjecture 13.}
The self-map $\psi_{\lambda}$
coincides with the  multiplication by  $p$ map on the  elliptic curve $\pi: \mathcal E_\lambda\to \mathbb P^1$ as the double cover ramifield on $S$ via $\pi$\\[.15cm]
We confirmed this conjecture via explicit computation for cases that the characteristic $p \leq 50$, see \cite{SYZ} Note that self-maps on $\mathbb P^1$ descended from $n$-multiplication by  maps on an elliptic curve $\pi: \mathcal E\to \mathbb P^1$
($n\in \mathbb Z$) are called Latt\`es maps (\cite{Mil}). If we vary $\lambda\in\mathbb P^1 \setminus \{0, 1, \infty\}$, we obtain a family of self-maps $\psi_{\lambda}$.  In the setting of complex dynamics, there is a rigidity theorem for stable families of self-maps  on $\mathbb P^1_{\mathbb C}$  due to McMullen \cite{Mc}. It states that any such family is either of Latt\`es type (i.e. constructed from $n$-multiplication map on $\mathcal E_\lambda$  by varying $\lambda$) or trivial (all its members are conjugate by M\"obius transformations).  We hope that the techniques in complex dynamics  will inspire us to investigate the self-map $\psi_\lambda$ induced from $p$-adic nonabelian Hodge theory by deforming the parameter $\lambda.$
\\[.1cm]
{\bf Problem 14.}   Study self maps arising from $p$-adic Higgs-de Rham flows in connection with the rigidiy theorem of McMullen on dynamical systems on the  complex projective line.
\\[.4cm]
{\bf Acknowledgment.} I would like to thank Steven Lu and Ruiran Sun for discussions on Riemann-Finsler
metric, Jinbang Yang on p-adic nonabelian Hodge theory and Junyi Xie on the works by Thurston and
Mcmullen on dynamic systems of self map on the projective line.\\[.4cm]
  {\bf  3. Bibliography.}\\


\begin{thebibliography}{KYZ20b}
  	
  	
  	
  	\bibitem[Abe13]{Ab}
  	Tomoyuki Abe.
  	\newblock Langlands correspondence for isocrystals and the existence of
  	crystalline companions for curves.
  	\newblock  J. Amer. Math. Soc. 31 (2018), 921-1057. 
  	
  	
  	\bibitem[Ara71]{A}
  	S.~Ju. Arakelov.
  	\newblock Families of algebraic curves with fixed degeneracies.
  	\newblock  Izv. Akad. Nauk SSSR Ser. Mat., 35:1269--1293, 1971.
  	
  	
  	
  	\bibitem[Ara19]{Ar}
  	Donu Arapura.
  	\newblock Kodaira-{S}aito vanishing via {H}iggs bundles in positive
  	characteristic.
  	\newblock  J. Reine Angew. Math., 755:293--312, 2019.
  	
  	\bibitem[BCS00]{CZ}
  	D.~Bao, S.-S. Chern, and Z.~Shen.
  	\newblock {\em An introduction to {R}iemann-{F}insler geometry}, volume 200 of
  	Graduate Texts in Mathematics.
  	\newblock Springer-Verlag, New York, 2000.
  	
  	\bibitem[CP15]{CP}
  	Fr\'{e}d\'{e}ric Campana and Mihai P\u{a}un.
  	\newblock Orbifold generic semi-positivity: an application to families of
  	canonically polarized manifolds.
  	\newblock  Ann. Inst. Fourier (Grenoble), 65(2):835--861, 2015.
  	
  	\bibitem[DF13]{DF}
  	Pierre Deligne and Yuval~Z. Flicker.
  	\newblock Counting local systems with principal unipotent local monodromy.
  	\newblock  Ann. of Math. (2), 178(3):921--982, 2013.
  	
  	\bibitem[{Den}18a]{De}
  	Ya~{Deng}.
  	\newblock Kobayashi hyperbolicity of moduli spaces of minimal
  	projective manifolds of general type (with the appendix by Dan Abramovich)
  	\newblock arXiv e-prints (2018) arXiv:1806.01666.
  	
  	
  	
  	\bibitem[DLSZ19]{DLSZ}
  	Ya Deng, Steven Lu, Ruiran Sun, and Kang Zuo.
  	\newblock Picard theorems for moduli spaces of polarized varieties, 2019.
  	
  	\bibitem[Dri81]{D}
  	V.~G. Drinfeld.
  	\newblock The number of two-dimensional irreducible representations of the
  	fundamental group of a curve over a finite field.
  	\newblock  Funktsional. Anal. i Prilozhen., 15(4):75--76, 1981.
  	\bibitem[Dri83]{D2}
  	V.G. Drinfed. 
  	\newblock Two-dimensional l-adic representations of the fundamental group of a curve over a finite
  	field and automorphic forms on GL(2).
  	\newblock  Amer. J. of Math. 105, No. 1 (1983), pp. 85–114.
  	
  	\bibitem[EG17]{EG}
  	H\'{e}l\`ene Esnault and Michael Groechenig.
  	\newblock Rigid connections and $f$-isocrystals, Acta Mathematica 225
  	1 (2020), 103–158.
  	
  	\bibitem[EK20]{EK}
  	H\'{e}l\`ene Esnault and  Moritz Kerz.
  	\newblock 	Density of arithmetic representations of function fields.
  	\newblock arXiv: 2005.12819.
  	
  	
  	
  	\bibitem[Fal83a]{F}
  	G.~Faltings.
  	\newblock Arakelov's theorem for abelian varieties.
  	\newblock  Invent. Math., 73(3):337--347, 1983.
  	
  	\bibitem[Fal83b]{F2}
  	G.~Faltings.
  	\newblock Endlichkeitss\"{a}tze f\"{u}r abelsche {V}ariet\"{a}ten \"{u}ber
  	{Z}ahlk\"{o}rpern.
  	\newblock  Invent. Math., 73(3):349--366, 1983.
  	
  	\bibitem[Fal89]{F3}
  	Gerd Faltings.
  	\newblock Crystalline cohomology and {$p$}-adic {G}alois-representations.
  	\newblock In Algebraic analysis, geometry, and number theory ({B}altimore,
  	{MD}, 1988), pages 25--80. Johns Hopkins Univ. Press, Baltimore, MD, 1989.
  	
  	\bibitem[Fal05]{F4}
  	Gerd Faltings.
  	\newblock A {$p$}-adic {S}impson correspondence.
  	\newblock Adv. Math., 198(2):847--862, 2005.
  	
  	\bibitem[FL82]{FL}
  	Jean-Marc Fontaine and Guy Laffaille.
  	\newblock Construction de repr\'esentations {$p$}-adiques.
  	\newblock  Ann. Sci. \'Ecole Norm. Sup. (4), 15(4):547--608 (1983), 1982.
  	
  	
  	\bibitem[G\"u60]{Gue}
  	Paul G\"{u}nther.
  	\newblock Einige {S}\"{a}tze \"{u}ber das {V}olumenelement eines {R}iemannschen
  	{R}aumes.
  	\newblock  Publ. Math. Debrecen, 7:78--93, 1960.
  	
  	\bibitem[GK73]{GK}
  	Phillip Griffiths and James King.
  	\newblock Nevanlinna theory and holomorphic mappings between algebraic
  	varieties.
  	\newblock  Acta Math., 130:145--220, 1973.
  	
  	\bibitem[Gri68a]{G1}
  	Phillip~A. Griffiths.
  	\newblock Periods of integrals on algebraic manifolds. {I}. {C}onstruction and
  	properties of the modular varieties.
  	\newblock  Amer. J. Math., 90:568--626, 1968.
  	
  	\bibitem[Gri68b]{G2}
  	Phillip~A. Griffiths.
  	\newblock Periods of integrals on algebraic manifolds. {II}. {L}ocal study of
  	the period mapping.
  	\newblock  Amer. J. Math., 90:805--865, 1968.
  	
  	\bibitem[Gri70]{G3}
  	Phillip~A. Griffiths.
  	\newblock Periods of integrals on algebraic manifolds. {III}. {S}ome global
  	differential-geometric properties of the period mapping.
  	\newblock  Inst. Hautes \'{E}tudes Sci. Publ. Math., (38):125--180, 1970.
  	
  	\bibitem[Hit87]	{H}
  	Nigel Hitchin
  	\newblock The Self-Duality Equations on a Riemann Surface.
  	\newblock Proceedings of the London Mathematical Society, 59-126, 1987.
  	
  	
  	
  	
  	\bibitem[JSZ20]{JSZ}
  	Ariyan Javanpeykar, Ruiran Sun, and Kang Zuo.
  	\newblock The Shafarevich conjecture revisited: Finiteness of pointed families
  	of polarized varieties, 2020.
  	
  	\bibitem[KL10]{KL}
  	S\'andor~J. Kov\'acs and Max Lieblich.
  	\newblock Boundedness of families of canonically polarized manifolds: a higher
  	dimensional analogue of {S}hafarevich's conjecture.
  	\newblock  Ann. of Math. (2), 172(3):1719--1748, 2010.
  	
  	\bibitem[Kon09]{Kon}
  	Maxim Kontsevich.
  	\newblock Notes on motives in finite characteristic.
  	\newblock In  Algebra, arithmetic, and geometry: in honor of {Y}u. {I}.
  	{M}anin. {V}ol. {II}, volume 270 of {\em Progr. Math.}, pages 213--247.
  	Birkh\"auser Boston, Inc., Boston, MA, 2009.
  	
  	\bibitem[KYZ20a]{KYZ2}
  	Raju Krishnamoorthy, Jinbang Yang, and Kang Zuo.
  	\newblock Finiteness of logarithmic crystalline representations, 2020.
  	
  	\bibitem[KYZ20b]{KYZ1}
  	Raju Krishnamoorthy, Jinbang Yang, and Kang Zuo.
  	\newblock A lefschetz theorem for crystalline representations, 2020.
  	
  	\bibitem[KYZ20c]{KYZ3}
  	Raju Krishnamoorthy, Jinbang Yang, and Kang Zuo.
  	\newblock Deformation theory of periodic Higgs-de Rham flows , 2020.
 
  	
  	
  	\bibitem[Lan15]{Lan}
  	Adrian Langer.
  	\newblock Bogomolov's inequality for {H}iggs sheaves in positive
  	characteristic.
  	\newblock  Invent. Math., 199(3):889--920, 2015.
  	
  	\bibitem[LSYZ19]{LSYZ}
  	Guitang Lan, Mao Sheng, Yanhong Yang, and Kang Zuo.
  	\newblock Uniformization of {$p$}-adic curves via {H}iggs--de {R}ham flows.
  	\newblock  J. Reine Angew. Math., 747:63--108, 2019.
  	
  	\bibitem[LSZ19]{LSZ1}
  	Guitang Lan, Mao Sheng, and Kang Zuo.
  	\newblock Semistable {H}iggs bundles, periodic {H}iggs bundles and
  	representations of algebraic fundamental groups.
  	\newblock  J. Eur. Math. Soc. (JEMS), 21(10):3053--3112, 2019.
  	
  	
  	
  	\bibitem[Mc87]{Mc}
  	Curt McMullen.
  	\newblock  Families of rational maps and iterative root-finding algorithms. 
  	\newblock 	Ann. of Math. (2) 125 (1987), no. 3, 467–493.
  	\bibitem[Mil68]{Mil68}
  	John Milnor.
  	\newblock A note on curvature and fundamental group.
  	\newblock  J. Differential Geometry, 2:1--7, 1968.
  	
  	\bibitem[Mil06]{Mil}
  	John Milnor.
  	\newblock On Latt\`es maps, Dynamics on the Riemann sphere.
  	\newblock Eur. Math. Soc., Z\"urich, 2006, pp. 9-43. 
  	
  	\bibitem[Moc96]{Mo}
  	Shinichi Mochizuki.
  	\newblock A theory of ordinary {$p$}-adic curves.
  	\newblock  Publ. Res. Inst. Math. Sci., 32(6):957--1152, 1996.
  	
  	\bibitem[OV07]{OV}
  	A. Ogus and V. Vologodsky.
  	\newblock Nonabelian {H}odge theory in characteristic {$p$}.
  	\newblock  Publ. Math. Inst. Hautes \'Etudes Sci., (106):1--138, 2007.
  	
  	\bibitem[PTW19]{PTW}
  	Mihnea Popa, Behrouz Taji, and Lei Wu.
  	\newblock Brody hyperbolicity of base spaces of certain families of varieties.
  	\newblock  Algebra Number Theory, 13(9):2205--2242, 2019.
  	
  	\bibitem[Sch13]{Sch}
  	Peter Scholze.
  	\newblock {$p$}-adic {H}odge theory for rigid-analytic varieties.
  	\newblock {\em Forum Math. Pi}, 1:e1, 77, 2013.
  	
  	\bibitem[Sim88]{S1}
  	Carlos~T. Simpson.
  	\newblock Constructing variations of {H}odge structure using {Y}ang-{M}ills
  	theory and applications to uniformization.
  	\newblock J. Amer. Math. Soc., 1(4):867--918, 1988.
  	
  	\bibitem[Sim10]{S2}
  	Carlos Simpson.
  	\newblock Iterated destabilizing modifications for vector bundles with
  	connection.
  	\newblock In  Vector bundles and complex geometry, volume 522 of 
  	Contemp. Math., pages 183--206. Amer. Math. Soc., Providence, RI, 2010.
  	
  	\bibitem[SYZ17]{SYZ}
  	Ruiran Sun, Jinbang Yang, and Kang Zuo.
  	\newblock Projective crystalline representations of \'etale fundamental groups
  	and twisted periodic higgs-de rham flow, 2017.
  	\newblock accepted by JEMS.
  	
  	\bibitem[VZ01]{VZ-0}
  	Eckart Viehweg and Kang Zuo.
  	\newblock On the isotriviality of families of projective manifolds over curves.
  	\newblock  J. Algebraic Geom., 10(4):781--799, 2001.
  	
  	\bibitem[VZ02]{VZ-2}
  	Eckart Viehweg and Kang Zuo.
  	\newblock Base spaces of non-isotrivial families of smooth minimal models.
  	\newblock In  Complex geometry ({G}\"ottingen, 2000), pages 279--328.
  	Springer, Berlin, 2002.
  	
  	\bibitem[VZ03]{VZ-1}
  	Eckart Viehweg and Kang Zuo.
  	\newblock On the {B}rody hyperbolicity of moduli spaces for canonically
  	polarized manifolds.
  	\newblock  Duke Math. J., 118(1):103--150, 2003.
  	
  	\bibitem[VZ05]{VZ-3}
  	Eckart Viehweg and Kang Zuo.
  	\newblock Complex multiplication, {G}riffiths-{Y}ukawa couplings, and rigidity
  	for families of hypersurfaces.
  	\newblock  J. Algebraic Geom., 14(3):481--528, 2005.
  	
  	\bibitem[WX07]{WX}
  	B.~Y. Wu and Y.~L. Xin.
  	\newblock Comparison theorems in {F}insler geometry and their applications.
  	\newblock  Math. Ann., 337(1):177--196, 2007.
  	
  	\bibitem[YZ20]{YZ}	
  	Jinbang Yang and Kang Zuo.
  	\newblock 	A note on p-adic Simpson correspondence.
  	\newblock ArXiv :2012.02058.
  	
  	
  	\bibitem[Zuo00]{Z}
  	Kang Zuo.
  	\newblock On the negativity of kernels of {K}odaira-{S}pencer maps on {H}odge
  	bundles and applications.
  	\newblock  Asian J. Math., 4(1):279--301, 2000.
  	\newblock Kodaira's issue.
  	
  \end{thebibliography}
\end{document}